\documentclass{article}
\usepackage[]{amsmath}
\usepackage[]{amsfonts}

\numberwithin{equation}{section}

\newtheorem{theorem}[equation]{Theorem}

\newtheorem{proposition}[equation]{Proposition}

\newtheorem{example}[equation]{Example}

\def\PPP{{\mathbb{P}}}

\begin{document}

\title{Pillow Degenerations of $K3$ Surfaces}
\author{C. Ciliberto, R. Miranda, and M. Teicher}
\maketitle

\section{Introduction}
In this article we construct a specific projective degeneration
of $K3$ surfaces of degree $2g-2$ in $\PPP^g$ to a union of
$2g-2$ planes, which meet in such a way that the combinatorics of
the configuration of planes is a triangulation of the $2$-sphere.
Abstractly, such degenerations are said to be Type III degenerations
of $K3$ surfaces, see \cite{kulikov}, \cite{pinkpers},
\cite{BGOD}.  Although the birational geometry of such degenerations
is fairly well understood, the study of projective degenerations
is not nearly as completely developed.

In \cite{clm}, projective
degenerations of $K3$ surfaces to unions of planes were constructed,
in which the general member was embedded by a primitive line bundle.
The application featured there was a computation of the rank of the
Wahl map for the general hyperplane section curve on the $K3$ surface.

In this article we construct degenerations for which the general member
is embedded by a \emph{multiple} of the primitive line bundle class.
The construction depends on two parameters, and we intend in follow-up
work to use these degenerations to compute braid monodromy for Galois
coverings, in the style of \cite{mt1} and \cite{mt2}.  We hope that
the freedom afforded by the additional discrete parameters in the
construction will yield interesting phenomena related to fundamental
groups.

The specific degenerations which we construct can be viewed as
two rectangular arrays of planes, joined along their boundary;
for this reason we have given them the name ``pillow'' degenerations.
They are described in Section \ref{pillow}.
Following that, in Section \ref{degenbranch}, we study the
degeneration of the general branch curve (for a general
projection of the surfaces to a plane) to a union of lines
(which is the ``branch curve'' for the union of planes).
In particular when the general branch curve is a plane curve
having only nodes and cusps as singularities, we describe
the degeneration of the nodes and the cusps to the configuration
of the union of lines.  This is critical information in the
application to the computation of the braid monodromy.

We are not aware of a modern reference for the statement
that the general branch curve for a linear projection of a surface
to a plane has only nodes and cusps as singularities.
In this article we will operate under the assumption that
this ``folklore'' statement is true and proceed.
The reader may wish to consult \cite{kk} for further information.
We have included a short section at the beginning of the article
deriving the characters of a general branch curve
(degree, number of nodes and cusps)
for the convenience of the reader,
under this assumption.

The authors are grateful to the NATO Scientific
Affairs Division, the Ministry of Science of Israel,
the Emmy Noether Research Institute of Mathematics
at Bar-Ilan University, and EAGER (the European
Union Research Network in Algebraic Geometry),
for financial support for the Workshop in Eilat at which
this article was completed.

\section{Characters of a General Branch Curve}
\label{branchcurve}

Here we briefly develop the formulas for the degree
and number of nodes and cusps on a general branch curve $B$
for a general projection of a smooth surface $S \subset \PPP^N$
to a general plane $\PPP^2$,
assuming that these are the only singularities.
These formulas are not new,
see for example \cite{enriques}, \cite{iversen},
but these standard references do much more, in either
outdated notation or with much more advanced techniques,
than are necessary for this more modest computation.  Hence
we thought it useful to include it here for completeness
and for the convenience of the reader.
The reader may also want to consult
\cite{kk}, \cite{mt1}, \cite{mtr}, and \cite{tr} for additional insight.

Denote by $\pi:S \to \PPP^2$ such a general projection.
Let $K$ and $H$ be the canonical and hyperplane classes of $S$
respectively.  Let $d$ be the degree of $S$ and $g(H)$ the
genus of a smooth hyperplane divisor.  The intersection numbers
$KH$ and $H^2$ are related to $d$ and $g$ by
\begin{equation}
d = H^2 \;\;\;\text{ and }\;\;\; 2g(H)-2 = H^2+KH.
\end{equation}
The degree of the finite map $\pi$ is equal to the degree $d$
of the surface $S$.
The degree $b$ of the branch curve may be easily computed by
noting that the pull-back of a line in $\PPP^2$ is a hyperplane
divisor; hence the Hurwitz formula gives
\[
2g(H)-2 = d(-2) + \deg(B)
\]
from which it follows that
\begin{equation}
b = \deg(B) = 2d + 2g(H) - 2 = 3d + KH.
\end{equation}

Let $R\subset S$ denote the ramification curve,
and denote by $R_0$ the residual curve (equal set-theoretically to
the closure of $\pi^{-1}(B) - R$).
$R$ is a smooth curve, and the mapping $\pi$, restricted to $R$,
is a desingularization of $B$.

Suppose that $B$ has $n$ nodes and $k$ cusps
and no other singularities.
Over a general smooth point of $B$, the map $\pi$ has $d-1$
preimages, one on the ramification curve and $d-2$ on the residual
curve.
Over each node of $B$, the ramification curve $R$ has two smooth
branches, and over each cusp, $R$ has one smooth branch.
Over a node of $B$, the residual curve $R_0$ meets $R$ once
transversally at each branch of $R$, and otherwise has $d-4$
nodes of its own.
Over a cusp of $B$, the residual curve $R_0$ meets $R$
twice at the point of $R$ lying over the cusp, and is smooth
there; it otherwise has $d-3$ cusps of its own.
In any case, over either a node or a cusp of $B$, there
are only $d-2$ preimages, instead of the $d-1$ preimages
over a general point of $B$.  Therefore, computing Euler numbers,
we see that
\begin{equation}
\label{e(RuR0)}
e(R\cup R_0) = e(\pi^{-1}(B)) = (d-1)e(B) -(n+k).
\end{equation}
The genus of the ramification curve $R$, being a
desingularization of the branch curve $B$, is
\[
g(R) = (b-1)(b-2)/2 - (n+k)
\]
using Pl\"ucker's formulas.
Its Euler number is therefore
\[
e(R) = 2-2g(R) = 2(n+k) - b^2 + 3b.
\]
Since $R$ and $B$ differ, topologically,
only over the nodes, we see that the Euler number of $B$
is
\[
e(B) = e(R) - n = n + 2k - b^2 + 3b.
\]
Letting $e(S)$ be the Euler number of the surface $S$,
we see that

$$e(S)= d[e(\PPP^2) - e(B)] + e(R\cup R_0)$$
$$= 3d -d e(B) + (d-1)e(B) -(n+k) \;\;\;\text{\rm using
\quad (\ref{e(RuR0)})}$$ $$= 3d -e(B) -n-k $$
$$= 3d -[n+2k-b^2+3b] - n - k $$
$$= 3d + b^2 - 3b -2n -3k,$$

so that
\begin{equation}
\label{2n+3k}
2n+3k = 3d + b^2 - 3b - e(S).
\end{equation}

Pulling back $2$-forms via $\pi$, we have the standard formula that
\[
K_S = \pi^*(K_{\PPP^2}) + R = -3H + R
\]
and since $bH = \pi^*(B)$, we see that
\[
2R+R_0 = \pi^*(B) = bH,
\]
so that, as classes on $S$,
\[
R = K + 3H \;\;\;\text{ and }\;\;\; R_0 = bH-2R = -2K +(b-6)H.
\]
Since $R$ and $R_0$ meet transversally at each of the two points
of $R$ over a node, and meet to order two at the point of $R$ lying
over a cusp, we see that $R\cdot R_0 = 2(n+k)$.
Therefore $2n + 2k = R\cdot R_0 = (K + 3H)(-2K + (b-6)H)$;
multiplying this out gives
\begin{equation}
\label{2n+2k}
2n + 2k = -2K^2 +(b-12)KH + (3b-18)H^2.
\end{equation}
Subtracting (\ref{2n+2k}) from (\ref{2n+3k}) gives
\[
k = 3d + b^2 - 3b - e(S) +2K^2 -(b-12)KH -(3b-18)d
\]
and then one can solve either expression for the number of nodes.
Simplifying the expressions somewhat leads to the following.

\begin{proposition}
Let $S$ be a smooth surface of degree $d$ in $\PPP^N$, and let
$\pi:S \to \PPP^2$ be a general projection.
Let $K$ and $H$ be the canonical and hyperplane classes of
$S$, respectively.
Let $B$ be the branch curve of the projection $\pi$,
which is assumed to be a plane curve of degree $b$
with $n$ nodes, $k$ cusps, and no other singularities.  Then:
\begin{enumerate}
\item[(a)] $\deg(\pi) = \deg(S) = d = H^2$.
\item[(b)] The degree of the branch curve $B$ is $b = 3d + KH$.
\item[(c)] The number of nodes of the branch curve $B$ is
\[
n = -3K^2 + e(S) + 24d + \frac{b^2}{2} -15b.
\]
\item[(d)] The number of cusps of the branch curve $B$ is
\[
k = 2K^2 - e(S) - 15d + 9b .
\]
\item[(e)] Under a general projection of the branch curve $B$
to a line, the number $t$ of turning points (simple branch points)
is
\[
t = e(S) - 3d + 2b.
\]
\end{enumerate}
\end{proposition}

The last computation of turning points is obtained from the
Hurwitz formula, applied to the ramification curve $R$,
noting that there are simple branch points for such a projection
at the points of $R$ lying over the cusps of $B$ also.

\begin{example}[Veronese Surfaces]
Let $S$ be the $r^{th}$ Veronese image of $\PPP^2$.
In this case, if $L$ denotes the line class of $S$, then
$L^2 = 1$, $K = -3L$, and $H = rL$;
hence $K^2 = 9$, $KH = -3r$, and $d = H^2 = r^2$.
The Euler number $e(S) = 3$.
Therefore
\[
\begin{array}{c}
b = 3r(r-1); \;\;\;
n = 3(r-1)(r-2)(3r^2+3r-8)/2; \\
k = 3(r-1)(4r-5); \;\;\;
t = 3(r-1)^2.
\end{array}
\]
\end{example}

\begin{example}[Rational Normal Scrolls]
Let $S$ be a rational normal scroll, e.g.
$\PPP^1\times\PPP^1$ embedded by the complete linear
system $H$ of type $(1,r)$.
The canonical class is of type $(-2,-2)$,
so that $K^2 = 8$, $KH = -2r-2$, and $d=H^2 = 2r$.
The Euler number $e(S)=4$. Therefore
\[
b = 4r-2; \;\;\;
n = 4(r-1)(2r-3); \;\;\;
k = 6r-6; \;\;\;
t = 2r.
\]
\end{example}

\begin{example}[Del Pezzo Surfaces]
Let $S$ be a Del Pezzo Surface of degree $d$ in $\PPP^d$,
for $3 \leq d \leq 9$.  Then $S$ is isomorphic to the
plane blown up at $9-d$ points; if $L$ denotes the class
of a line, and $E$ the sum of the classes of the $9-d$
exceptional divisors, then
$L^2 = 1$, $LE = 0$, and $E^2 = d-9$;
also $K = -3L+E$, and $H = -K$, so that
$K^2 = H^2 = d$, and $KH = -d$.
The Euler number $e(S) = 12-d$.  Therefore
\[
b = 2d; \;\;\;
n = 2(d-2)(d-3) =2d^2 - 10d +12; \;\;\;
k = 6(d-2); \;\;\;
t = 12.
\]
\end{example}

\begin{example}[$K3$ Surfaces]
Let $S$ be a $K3$ surface of degree $d = 2g-2$
in $\PPP^g$.  The canonical class is trivial,
so that $K^2=KH=0$.
The Euler number $e(S) = 24$.  Therefore
\[
\begin{array}{c}
b = 6g-6; \;\;\;
n = 6(g-2)(3g-7) = 18g^2-78g+84; \\
k = 24(g-2); \;\;\;
t = 6g+18.
\end{array}
\]
\end{example}

\section{Construction of the Pillow Degeneration}
\label{pillow}

A non-hyperelliptic $K3$ surface of genus $g \geq 3$
can be embedded by the sections of a very ample line bundle
as a smooth surface of degree $2g-2$ in $\PPP^g$.
When the line bundle generates the Picard group of the $K3$ surface,
the embedded $K3$ surface can be degenerated to a union of $2g-2$ planes
in a variety of ways (see for example \cite{clm}).
In this section we will describe a degeneration,
which we call the \emph{pillow} degeneration,
which smooths to a $K3$ surface whose Picard group
is generated by a sub-multiple of the hyperplane class.

Fix two integers $a$ and $b$ at least two;
set $g = 2ab+1$.
The number of planes in the pillow degeneration
is then $2g-2 = 4ab$.

This projective space has $g+1 = 2ab+2$ coordinate points,
and each of the $4ab$ planes is obtained as the span of three
of these.  The sets of three are indicated in Figure 1,
which describes the bottom part of the ''pillow''
and the top part of the ''pillow'',
which are identified along the boundaries of the two
configurations.
The reader will see that the boundary is a cycle of $2a+2b$ lines.

\begin{figure}[h]
\caption{Configuration of Planes, Top and Bottom}
\begin{picture}(500,180)
\put(70,170){Top}
\put(0,160){Boundary Points Labeled From}
\put(0,150){1 through 2a+2b, clockwise;}
\put(0,140){Interior Points Labeled from}
\put(0,130){2a+2b+1 through ab+a+b+1}
\put(27,105){1}
\put(37,105){2}
\put(52,105){$\cdots$}
\put(117,105){a}
\put(132,105){a+1}
\put(0,90){2a+2b}\put(132,90){a+2}
\put(132,80){a+3}
\put(10,60){$\vdots$}\put(132,60){$\vdots$}
\put(0,10){a+2b+1}\put(127,10){a+b+1}
\multiput(30,20)(10,0){10}{\frame{\line(1,1){10}}}
\multiput(30,30)(10,0){10}{\frame{\line(1,1){10}}}
\multiput(30,40)(10,0){10}{\frame{\line(1,1){10}}}
\multiput(30,50)(10,0){10}{\frame{\line(1,1){10}}}
\multiput(30,60)(10,0){10}{\frame{\line(1,1){10}}}
\multiput(30,70)(10,0){10}{\frame{\line(1,1){10}}}
\multiput(30,80)(10,0){10}{\frame{\line(1,1){10}}}
\multiput(30,90)(10,0){10}{\frame{\line(1,1){10}}}
\put(250,170){Bottom}
\put(190,160){Boundary Points Labeled From}
\put(190,150){1 through 2a+2b, clockwise;}
\put(190,140){Interior Points Labeled from}
\put(190,130){ab+a+b+2 through 2ab+2}
\put(207,105){1}
\put(217,105){2}
\put(232,105){$\cdots$}
\put(297,105){a}
\put(312,105){a+1}
\put(180,90){2a+2b}\put(312,90){a+2}
\put(312,80){a+3}
\put(190,60){$\vdots$}\put(312,60){$\vdots$}
\put(180,10){a+2b+1}\put(307,10){a+b+1}
\multiput(210,20)(10,0){10}{\frame{\line(1,-1){10}}}
\multiput(210,30)(10,0){10}{\frame{\line(1,-1){10}}}
\multiput(210,40)(10,0){10}{\frame{\line(1,-1){10}}}
\multiput(210,50)(10,0){10}{\frame{\line(1,-1){10}}}
\multiput(210,60)(10,0){10}{\frame{\line(1,-1){10}}}
\multiput(210,70)(10,0){10}{\frame{\line(1,-1){10}}}
\multiput(210,80)(10,0){10}{\frame{\line(1,-1){10}}}
\multiput(210,90)(10,0){10}{\frame{\line(1,-1){10}}}
\end{picture}
\end{figure}
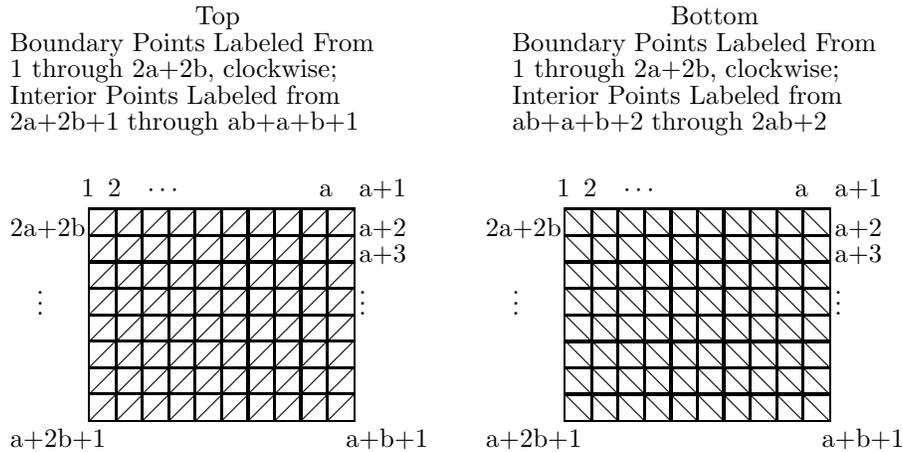

Note that no three of the planes meet in a line.
Also note that the set of bottom planes
lies in a projective space of dimension $ab+a+b$,
as do the set of top planes;
these two projective spaces meet
exactly along the span of the $2a+2b$ boundary points,
which has dimension $2a+2b-1$.
Finally note that the four corner points of the pillow
degeneration
(labeled $1$, $a+1$, $a+b+1$, and $a+2b+1$)
are each contained in three distinct planes,
while all other points are each contained in six planes.
This property, that the number of lines and planes
incident on each of the points is bounded,
is important for the later computations,
and is a feature of the pillow degeneration
that is not available in other previous degenerations.

We will call such a configuration of planes a 
\emph{pillow of bidegree $(a,b)$}.

\begin{theorem}
For any $a$ and $b$ at least $2$, the pillow of bidgree $(a,b)$
is a degeneration of a smooth $K3$ surface
of degree $4ab$ in a projective space
of dimension $g = 2ab+1$.
If $c = g.c.d(a,b)$, then the general such $K3$ surface
will have Picard group generated by a line bundle $L$
such that $cL$ is the hyperplane bundle.
\end{theorem}

\medskip
The proof of the Theorem will be made
in three steps.
First we will exhibit a degeneration of the $K3$
surface to a union of two rational surfaces,
each isomorphic to $\PPP^1\times\PPP^1$,
embedded via the sections of the linear system
of bidegree $(a,b)$.
These two rational surfaces will meet along an
elliptic normal curve which is anticanonical in each.
Secondly we will simultaneously degenerate
each rational surface
to a union of $ab$ quadrics,
resulting in a total of $2ab$ quadrics.
Finally we will degenerate each quadric to a union of two planes.

\medskip\noindent
{\bf Proof:} (Step One:)
Note that the sections of the linear system of bidegree $(a,b)$
embed $\PPP^1\times\PPP^1$ as a surface
in a projective space of dimension $ab+a+b$.
Choose an anticanonical divisor (of bidegree $(2,2)$)
which is a smooth elliptic curve;
it is mapped by the above embedding
to an elliptic normal curve in a subspace
of dimension $2a+2b-1$.

In our original ambient space of dimension $2ab+1$,
choose two subspaces of dimension $ab+a+b$
which meet along a subspace of dimension $2a+2b-1$.
Make the above identical construction of the $\PPP^1\times\PPP^1$
in each of the two subspaces, taking care
to have the two elliptic normal curves identified
in the intersection subspace.

This union $R$ of the two rational surfaces
is a degeneration of an embedded $K3$ surface,
by an argument identical to that presented for
Theorems 1 and 2 of \cite{clm},
which we will not repeat in detail here.
Briefly, one first checks via standard calculations
that $H^0(N_R)$ has dimension $g^2+2g+19$
and that $H^1(N_R) = H^2(N_R) = 0$.
Secondly, the natural map from $H^0(N_R)$ to $H^0(T^1)$
is seen to be surjective.
This is sufficient to prove that $R$ represents a smooth
point of its Hilbert scheme, whose general member is a smooth
$K3$ surface of degree $2g-2$.

\medskip\noindent
(Step Two:)
The second step can be achieved as in \cite{mt1} by observing that
each $\PPP^1\times\PPP^1$ can be degenerated to a union of $ab$
quadrics by degenerating
the first coordinate $\PPP^1$ to a chain of $a$ lines,
and the second coordinate $\PPP^1$ to a chain of $b$ lines.
In this degeneration the elliptic curve
degenerates to a cycle of $2(a+b)$ lines.
This degeneration is made simultaneously for each of the
two $\PPP^1\times\PPP^1$'s, resulting in a degeneration
to a union of $2ab$ quadrics.
This configuration of $2ab$ quadrics meet as in Figure 1,
without the diagonal lines: if one removes the diagonal
lines from Figure 1 we obtain $2ab$ rectangles,
each indicating a quadric.  Each of these quadrics
meets the others along a cycle of four lines
(two vertical and two horizontal).

\medskip\noindent
(Step Three:) Finally degenerate each quadric to a union of
two planes, as in Figure 1.  These degenerations can
be executed completely independently of course,
and it is elementary to see that this can be done
keeping the four lines along which any one of the quadrics
meet the others fixed.

\medskip\noindent
(Step Four:)
Finally note that if $c \neq 1$,
the pillow degeneration of bidegree $(a,b)$
is a degeneration of the $c$-uple embedding of the pillow
degeneration of bidegree $(a/c,b/c)$.
To see this, one uses the standard triangular degeneration
of the Veronese embedding of the plane as described in \cite{mt2}.

The final point to check is that the general $K3$
surface in this $19$-dimensional family has Picard group
generated by $(1/c)H$, where $H$ is the hyperplane class.
Since we have a $19$-dimensional
family of $K3$ surfaces, the only question to be decided is
which sub-multiple of the hyperplane system is the generator
of the general Picard group.
The maximum possible is the g.c.d $c$.
Since the pillow is a $c$-uple Veronese,
the hyperplane class is at least a $c$-fold
multiple, and since it cannot be any more,
this shows that the Picard group is generated by $(1/c)H$.
This completes the proof of the Theorem.

\rightline{Q.E.D}

\medskip
Note that in this degeneration,
the horizontal and vertical lines appear first,
and the diagonal lines appear second.

\section{The Degeneration of the Branch Curve}
\label{degenbranch}

We assume that we are in a general enough
situation that 
for a generic projection
of a K3 surface of degree $g$ in $\PPP^g$ to a plane,
the branch curve is a curve of degree $6g-6$,
having $6(g-2)(3g-7)$ nodes and $24(g-2)$ cusps
and no other singularities;
these numbers were presented in Section \ref{branchcurve}.
If one projects this branch curve onto a general line,
the projection will have $6g+18$ simple branch points.

It is our goal in this section to describe
how these nodes, cusps, and branch points degenerate
in a pillow degeneration.

Firstly, since the pillow degeneration consists entirely of planes,
under a general projection each plane will map isomorphically onto
the target plane.  Therefore the degenerate branch curve is composed of
the $3g-3$ planar lines which are the images of the $3g-3$
double lines of the pillow degeneration where two planes meet.
Each of the $3g-3$ planar lines have multiplicity two in the
limit branch curve.

We see therefore that the general branch curve (of degree $6g-6$)
degenerates as a curve to the $3g-3$ planar lines, each doubled.
Our next task is to describe the degeneration of the
nodes, cusps and branch points of the general branch curve.
In any case it is clear that these distinguished points 
of the general branch curve
can only go to points of the $3g-3$ planar lines.

Secondly, it is elementary to compute that there are
$(9/2)g^2 - (51/2) g + 39$ pairs of disjoint lines
in the pillow degeneration.
Each of these pairs of disjoint lines gives rise to an
intersection of two planar line components of the limit branch curve.
We refer to these points as \emph{$2$-points} of the configuration
of the $3g-3$ planar lines.

In addition to these $2$-points,
we have exactly four \emph{$3$-points},
corresponding to the projection of the four points in
the pillow degeneration where $3$ planes (and $3$ double lines)
meet.
Finally we have $g-3$ \emph{$6$-points}
corresponding to the projection of the $g-3$ points in
the pillow degeneration where $6$ planes (and $6$ double lines)
meet.
At any one of these $n$-points ($n=2$, $3$, or $6$)
exactly $n$ of the $3g-3$ planar lines meet;
moreover at no other point of the plane do any of these lines meet.

In the degeneration of the general branch curve
to this configuration of $3g-3$ double lines,
each of the nodes, cusps, and branch points can degenerate
either to a $2$-point, a $3$-point, a $6$-point,
or a smooth point of one of the $3g-3$ lines.
With the above terminology, we can now describe how many
nodes, cusps, and branch points degenerate to each of these types of
points.

\begin{theorem}
In the pillow degeneration of a $K3$ surface of degree $g$ in $\PPP^g$,
the nodes, cusps, and branch points of the general branch curve
degenerate to the $2$-points, $3$-points, $6$-points,
and other smooth points of lines
according to the following table:

\begin{tabular}{c|c|c|c|c}
Object & Number & Branch & Nodes & Cusps \\
Type & & Points & & \\ \hline
Lines & $3g-3$ & $0$ & $0$ & $0$ \\ \hline
$3$-points & $4$ & $9$ & $0$ & $6$ \\ \hline
$6$-points & $g-3$ & $6$ & $24$ & $24$ \\ \hline
$2$-points & $\frac{9}{2}g^2-\frac{51}{2}g+39$ & $0$ & $4$ & $0$ \\ \hline
Totals: & & $6g+18$ & $18g^2-78g+84$ & $24(g-2)$
\end{tabular}

In particular no node, cusp, or branch point degenerates
to a smooth point of any of the $3g-3$ double lines of the
limit branch curve.
\end{theorem}

\medskip\noindent{\bf Proof:}
We first look at the row of the table for the $2$-points.
Since each of the planar lines have multiplicity two in the branch curve,
this crossing point actually is a limit of $4$ nodes of the general branch
curve
(the $4$ nodes appearing as the four intersection points of two pairs
of lines).
No cusp or branch point of the general branch curve
has this crossing point as a limit in general,
since these points are created by the projection of unrelated disjoint
lines
in the union of planes in $\PPP^g$.

We now turn our attention to the images of the multiple points
of the pillow degeneration where $n$ planes (and $n$ double lines)
meet at one point.  We assume that $3 \leq n \leq 6$ in what follows.
(In the pillow degeneration we have $n=3$ or $n=6$ only.)
Under the generic planar projection, such points
go to intersections of $n$ of the corresponding planar lines.
We will refer to these as \emph{$n$-points} of the limit branch curve.

In order to analyze the number of nodes, cusps, and branch points
of the general curve which go to these $n$-points,
we make a local analysis near the multiple point of the union of planes.
There are $n$ planes incident to this multiple point,
and they together span a $\PPP^n$.
Locally this collection of $n$ planes in $\PPP^n$
smooths to a Del Pezzo surface of degree $n$.
In a generic projection for such a Del Pezzo, 
the branch curve has degree $2n$,
with $2(n-2)(n-3)$ nodes and $6n-12$ cusps;
the number of simple branch points for this curve
under generic projection to a line is $12$.

The limit branch curve corresponding to the degeneration
of the Del Pezzo to the union of $n$ planes is a union of
$n$ lines concurrent at a point $p$, the images of the $n$ lines through
the multiple point.

A partial smoothing of the union of $n$ planes
may be obtained by taking two adjacent planes and smoothing
them to a quadric surface.
The corresponding smoothing of the limit branch curve
smooths exactly one of the $n$ lines to a conic,
which is necessarily tangent to two adjacent lines.
As the conic degenerates to the (double) line $L$,
we see that no nodes of the general branch curve go to
any point of $L$ which is not $p$, and no cusps do either.
The conic has two general branch points for a projection to a line,
and one of these branch points goes to $p$ and one does not.

This local analysis of this partial smoothing shows that
in a complete smoothing to the Del Pezzo,
no node can go to a point of any line except the concurrent point $p$,
and neither can any cusp.
Therefore all of the $2(n-2)(n-3)$ nodes degenerate to the concurrent
point,
and all of the $6n-12$ cusps do too.
Moreover, of the $12$ branch points for the general curve,
all but $n$ of them go to the concurrent point $p$.
(The other $n$ go to one on each line.)

In the cases $n=3$ and $n=6$ of interest in the pillow degeneration,
the above analysis shows that arbitrarily close to a $3$-point
there are $9 = 12-3$ branch points, and no nodes and $6$ cusps.
Arbitrarily close to a $6$-point there are
$6=12-6$ branch points, and $24$ nodes and $24$ cusps.
This gives the entries in the $3$-point and $6$-point
rows of the table.

If we now total the number of branch points,
nodes and cusps which degenerate to these multiple points,
we obtain the values in the last row of the table.
Since these are exactly the number of branch points,
nodes, and cusps of the general curve,
we must have accounted for all of the branch points, nodes, and cusps
already.  In particular there are none left to degenerate to smooth
points of the double lines.

This completes the proof of the Theorem.

\rightline{Q.E.D}

\end{document}